\documentclass[12pt]{article}
\usepackage{amsmath,amsfonts,amsthm,amssymb}
\newcommand{\F}{\mathbb F_q}
\newcommand{\K}{\overline{K}_c}
\newcommand{\AG}{\mathfrak A}

\newcommand{\FF}{\mathcal F}
\newcommand{\D}{\mathcal D}
\newcommand{\KT}{\overline{K}_c\{ t\}}
\DeclareMathOperator{\gr}{gr}
\DeclareMathOperator{\dist}{dist}

\numberwithin{equation}{section}
\begin{document}
\newtheorem{teo}{Theorem}[section]
\pagestyle{plain}
\title{Hypergeometric Functions and Carlitz Differential Equations
over Function Fields}
\author{Anatoly N. Kochubei\\ 
\footnotesize Institute of Mathematics,\\ 
\footnotesize National Academy of Sciences of Ukraine,\\ 
\footnotesize Tereshchenkivska 3, Kiev, 01601 Ukraine
\\ \footnotesize E-mail: \ kochubei@i.com.ua}
\date{}
\maketitle
\newpage
\begin{abstract}
The paper is a survey of recent results in analysis of additive 
functions over function fields motivated by applications to 
various classes of special functions including Thakur's 
hypergeometric function. We consider basic notions and results of 
calculus, analytic theory of differential equations with Carlitz 
derivatives (including a counterpart of regular singularity), 
umbral calculus, holonomic modules over the Weyl-Carlitz ring.
\end{abstract}
\vspace{2cm}
\newpage

\section{INTRODUCTION}

Let $K$ be the field of formal Laurent series
$t=\sum\limits_{j=N}^\infty \xi_jx^j$ with coefficients $\xi_j$ 
from the Galois field $\F$, $\xi_N\ne 0$ if $t\ne 0$, $q=p^\upsilon $, 
$\upsilon \in \mathbf Z_+$,
where $p$ is a prime number. It is well known that any non-discrete 
locally compact field of characteristic $p$ is isomorphic to such $K$.
The absolute value on $K$ is given by $|t|=q^{-N}$, $|0|=0$. The ring of
integers $O=\{ t\in K:\ |t|\le 1\}$ is compact in the topology
corresponding to the metric $\dist (t,s)=|t-s|$.  The
absolute value $|\cdot |$ can be extended in a unique way onto the 
completion $\K$ of an algebraic closure of $K$.

Analysis over $K$ and $\K$ initiated  
by Carlitz \cite{C1} and developed subsequently by Wagner, Goss, 
Thakur, the author, and others (see the bibliography in 
\cite{G2,Th3}) is very different from the classical calculus. The 
new features begin with an appropriate version of the factorial 
invented by Carlitz -- since the usual factorial $i!$, seen as an 
element of $K$, vanishes for $i\ge p$, Carlitz introduced the new 
one as
\begin{equation} 
D_i=[i][i-1]^q\ldots [1]^{q^{i-1}},\quad [i]=x^{q^i}-x\ (i\ge 1),\ 
D_0=1.
\end{equation}

An important feature is the availability of many non-trivial 
$\F$-linear functions, that is such functions $f$ defined on 
$\F$-subspaces $K_0\subset K$ that
$$
f(t_1+t_2)=f(t_1)+f(t_2),\quad f(\alpha t)=\alpha f(t),
$$
for any $t,t_1,t_2\in K_0$, $\alpha \in \F$. Such are, for 
example, polynomials and power series of the form $\sum 
a_kt^{q^k}$, in particular, the Carlitz exponential
\begin{equation}
e_C(t)=\sum\limits_{n=0}^\infty \frac{t^{q^n}}{D_n},\quad |t|<1,
\end{equation} 
and its composition inverse, the Carlitz logarithm
\begin{equation}
\log_C(t)=\sum\limits_{n=0}^\infty (-1)^n\frac{t^{q^n}}{L_n},\quad |t|<1,
\end{equation}
where $L_n=[n][n-1]\cdots [1]$ ($n\ge 1$), $L_0=1$. The notion of 
the Carlitz exponential obtained a wide generalization in the 
theory of Drinfeld modules (see \cite{G2,Th3}). On the other hand, 
in various problems going beyond the class of $\F$-linear 
functions, an extended version of the Carlitz factorial (and its 
Gamma function interpolations) is used, so that $D_n$ can be seen 
as ``an $\F$-linear part'' of the full factorial; see 
\cite{G2,Th3} and references therein for the details.

Among other special classes of $\F$-linear functions there are 
various polynomial systems (see below), an analog of the Bessel 
functions \cite{C2,SL}, and Thakur's hypergeometric function 
\cite{Th1,Th2,Th3}. The latter is defined as follows.

For $n\in \mathbb Z_+$, $a\in \mathbb Z$, denote
\begin{equation}
(a)_n=\begin{cases}
D_{n+a-1}^{q^{-(a-1)}}, & \text{if $a\ge 1$;}\\
L_{-a-n}^{-q^n}, & \text{if $a\le 0,n\le -a$;}\\
0, & \text{if $a\le 0,n>-a$}.\end{cases}
\end{equation} 
Then, for $a_i,b_i\in \mathbb Z$, such that the series below makes 
sense, we set
\begin{equation}
{}_rF_s(a_1,\ldots ,a_r;b_1,\ldots 
,b_s;z)=\sum\limits_{n=0}^\infty \frac{(a_1)_n\cdots 
(a_r)_n}{(b_1)_n\cdots (b_s)_nD_n}z^{q^n}.
\end{equation} 

Thakur \cite{Th1,Th2,Th3} has carried out a thorough investigation 
of the functions (1.5) and obtained analogs of many properties 
known for the classical situation. In particular, he found an 
analog of the hypergeometric differential equation. Its main 
ingredients are the difference operator
$$
(\Delta u)(t)=u(xt)-xu(t)
$$
(an inner derivation of composition rings of $\F$-linear 
polynomials or more general $\F$-linear functions)
introduced by Carlitz \cite{C1}, the nonlinear ($\F$-linear) operator 
$d=\sqrt[q]{}\circ \Delta$, and the $\F$-linear Frobenius operator
$\tau u=u^q$. For example, the function $y={}_2F_1(a,b;c;z)$ is a 
solution of the equation
\begin{equation}
(\Delta -[-a])(\Delta -[-b])y=d(\Delta -[1-c])y.
\end{equation} 
Here we touch only a part of Thakur's results (he considered also 
hypergeometric functions corresponding to other places of $\F 
(x)$, a version of (1.5) with parameters from $K$ and its 
extensions etc).

The Carlitz exponential $e_C$ satisfies a much simpler equation of 
the same kind:
\begin{equation}
de_C=e_C,
\end{equation}
so that the operator $d$ may be seen as an analog of the 
derivative. The operator $\tau$ is an analog of the 
multiplication by $t$, so that $\Delta$ is the counterpart of 
$t\dfrac{d}{dt}$.

The same operators appear in the positive characteristic analogs 
of the canonical commutation relations of quantum mechanics 
\cite{K98,K99}. In the analog of the Schr\"odinger representation 
we consider, on the Banach space $C_0(O,\K )$ of continuous 
$\F$-linear functions on $O$, with values from $\K$ (with the 
supremum norm), the ``creation and annihilation operators''
$$
a^+=\tau -I,\quad a^-=d
$$
($I$ is the identity operator). Then
\begin{equation}
a^-a^+-a^+a^-=[1]^{1/q}I,
\end{equation}
the operator $a^+a^-$ possesses the orthonormal (in the 
non-Archimedean sense \cite{Sch}) eigenbasis $\{ f_i\}$,
\begin{equation}
(a^+a^-)f_i=[i]f_i,\quad i=0,1,2,\ldots ;
\end{equation}
$a^+$ and $a^-$ act upon the basis as follows:
\begin{equation}
a^+f_{i-1}=[i]f_i,\ \ a^-f_i=f_{i-1},\ i\ge 1;\ a^-f_0=0.
\end{equation}

Here $\{ f_i\}$ is the sequence of normalized Carlitz polynomials
\begin{equation}
f_i(s)=D_i^{-1}\prod \limits _{\genfrac{}{}{0pt}{1}{m\in 
\F [x]}{\deg m<i}}(s-m) \quad (i\ge 1),\quad f_0(s)=s,
\end{equation}
which forms an orthonormal basis in $C_0(O,\K )$. The spectrum of 
the ``number operator'' $a^+a^-$ is the set of elements $[i]$, so 
that even this notation (proposed by Carlitz in 1935) becomes 
parallel to the usual quantum mechanical situation.

An analog of the Bargmann-Fock representation is obtained if we 
consider the operators of almost the same form,
$$
\tilde{a}^+=\tau,\quad \tilde{a}^-=d,
$$
but on the Banach space $H$ of power series 
$u(t)=\sum\limits_{n=0}^\infty a_n\frac{t^{q^n}}{D_n}$ with 
$a_n\in \K$, $a_n\to 0$ as $n\to \infty$. These new operators 
satisfy the same relations (1.8)-(1.10), but this time instead of 
the Carlitz polynomials $f_n$ we get the eigenfunctions 
$\tilde{f}_n=\frac{t^{q^n}}{D_n}$.

The above results motivated the author to begin to develop 
analysis and theory of differential equations for $\F$-linear 
functions over $K$ and $\K$, that is for the case which can be 
seen as a concentrated expression of features specific for the 
analysis in positive characteristic. This paper is a brief survey 
of some achievements in this direction. In particular, we consider 
the counterparts of the basic notions of calculus, analytic theory 
of differential equations (in the regular case and the case of 
regular singularity), their applications to some special 
functions, like the power function, logarithm and polylogarithms, Thakur's 
hypergeometric function etc. An umbral calculus and a theory of 
holonomic modules are initiated for this case. Like in the 
classical situation (see \cite{Cart}), it is shown that some basic 
objects of the function field arithmetic generate holonomic modules.

Note that some of the results can be easily extended to the case 
where the base field is a completion of $\F (x)$ with respect to a 
finite place determined by an irreducible polynomial $\pi \in \F 
[x]$ (the field $K$ corresponds to $\pi (x)=x$); for some details 
see \cite{Kzeta}. The situation is different for the ``infinite'' 
place widely used in function field arithmetic (see \cite{Th3}). 
In this case some of the basic objects behave in a quite different 
way -- absolute values of the Carlitz factorials $D_n$ grow, as 
$n\to \infty$, the Carlitz exponential is an entire function, the 
Carlitz polynomials do not form an orthonormal basis etc. A 
thorough investigation of properties of the Carlitz differential 
equations for this situation has not been carried out so far.

\section{Calculus}

{\bf 2.1.} Higher Carlitz operators $\Delta^{(n)}$ are introduced 
recursively,
\begin{equation}
\left( \Delta ^{(n)}u\right) (t)=\Delta ^{(n-1)}u(xt)-
x^{q^{n-1}}\Delta ^{(n-1)}u(t),\quad n\ge 2.
\end{equation}
For $n=1$, the formula (2.1) coincides with the definition of 
$\Delta =\Delta^{(1)}$, if we set $\Delta^{(0)}=I$.

The first application of these operators is the reconstruction 
formula \cite{K99} for the coefficients $a_n$ of a power series 
$u\in H$. Note that the classical formula does not make sense here 
because it contains the expression $u^{(n)}(t)/n!$ where both the 
numerator and denominator vanish.

\medskip
\begin{teo}
If $u\in H$, then
$$
a_n=\lim \limits _{t\to 0}\frac{\Delta ^{(n)}u(t)}{t^{q^n}},\quad
n=0,1,2,\ldots .
$$
\end{teo}

\medskip
For a continuous non-holomorphic $\F$-linear function $u$ the 
behaviour of the functions
$$
\mathfrak D^ku(t)=t^{-q^k}\Delta ^{(k)}u(t),\quad t\in O\setminus
\{ 0\} ,
$$
near the origin measures the smoothness of $u$. We say that 
$u\in C_0^{k+1}(O,\K )$ if $\mathfrak D^ku$ can
be extended to a continuous function on $O$. This includes the 
case ($k=0$) of differentiable functions.

The next theorem proved in \cite{K99} gives a characterization of 
the above smoothness in terms of coefficients of the 
Fourier-Carlitz expansion. It includes, as a particular case 
($k=0$), the characterization of differentiable $\F$-linear functions 
obtained by Wagner \cite{W}.

\medskip
\begin{teo}
A function $u=\sum \limits _{n=0}^\infty c_nf_n\in C_0(O,\K )$
belongs to $C_0^{k+1}(O,\K )$ if and only if
$$
q^{nq^k}|c_n|\to 0\quad \mbox{for }n\to \infty .
$$
In this case
$$
\sup \limits _{t\in O}|\mathfrak D^ku(t)|=\sup \limits 
_{n\ge k}q^{(n-k)q^k}|c_n|.
$$
\end{teo}

\medskip
For a generalization to some classes of not necessarily 
$\F$-linear functions see \cite{Y2}.

Similarly \cite{K99}, a function $u$ is analytic  on the ball $O$
(that is, $u(t)=\sum a_it^{q^i}$, $a_i\to 0$) if and only if 
$q^{\frac{q^n}{q-1}}|c_n|\to 0$, as $n\to \infty$. A more refined 
result by Yang \cite{Y1}, useful in many applications, which makes 
it possible to find an exact domain of analyticity, is as follows 
(again we consider only $\F$-linear functions while in \cite{Y1} a 
more general class is studied).

\medskip
\begin{teo}[Yang]
A function $u=\sum \limits _{n=0}^\infty c_nf_n\in C_0(O,\K )$ is 
locally analytic if and only if
\begin{equation}
\gamma =\liminf\limits_{n\to \infty }\left\{ -q^{-
n}\log_q|c_n|\right\} >0,
\end{equation}
and if (2.2) holds, then $u$ is analytic on any ball
of the radius $q^{-l}$, 
$$
l=\max (0,[-(\log (q-1)+\log \gamma )/\log q]+1).
$$
\end{teo}

\medskip
{\bf 2.2.} Viewing $d$ as a kind of a derivative, it is natural to 
introduce an antiderivative $S$ setting $Sf=u$ where $u$ is a 
solution of the equation $du=f$, with the normalization $u(1)=0$. 
It is easy to find $Sf$ explicitly if $f$ is given by its 
Fourier-Carlitz expansion (see \cite{K99}).

Next, we introduce a Volkenborn-type integral of a function
$f\in C_0^1(O,\K )$ (see \cite{Sch} for a similar integration 
theory over $\mathbb Z_p$) setting
$$
\int \limits _{O}f(t)\,dt\stackrel{\mbox{{\footnotesize def}}}{=}\lim \limits
_{n\to \infty }\frac{Sf(x^n)}{x^n}=(Sf)'(0).
$$

The integral is a $\F$-linear continuous
functional on $C_0^1(O,\K )$,
$$
\int \limits _{O}cf(t)\,dt=c^q\int \limits _{O}f(t)\,dt,\quad c\in \K ,
$$
possessing the following ``invariance'' property
(related, in contrast to the case of $\mathbf Z_p$, to the
multiplicative structure):
$$
\int \limits _{O}f(xt)\,dt=x\int \limits _{O}f(t)\,dt-f^q(1).
$$

Calculating the integrals of some important functions we obtain 
new relations between them. In addition to the Carlitz exponential 
$e_C$ and the Carlitz polynomials $f_n$ (see (1.2) and (1.11)), we mention
the Carlitz module function
\begin{equation}
C_s(z)=\sum \limits _{i=0}^\infty f_i(s)z^{q^i},\quad s\in O,|z|<1.
\end{equation}
Note that if $s\in \F [x]$, then only the terms with $i\le \deg s$ 
are different from zero in (2.3).

We have
$$
\int \limits _{O}t^{q^n}\,dt=-\frac{1}{[n+1]}, \quad n=0,1,2,\ldots ;
$$
$$
\int \limits _{O}f_n(t)\,dt=\frac{(-1)^{n+1}}{L_{n+1}},\quad n=0,1,2,\ldots ;
$$
$$
\int \limits _{O}C_s(z)\,ds=\log _C(z)-z,\quad z\in K,\ |z|<1;
$$
$$
\int \limits _Oe_C(st)\,ds=t-e_C(t),\quad t\in K,|t|<1.
$$
For the proofs see \cite{K99}.

\section{Differential Equations for $\F$-Linear Functions}

{\bf 3.1.} Let us consider function field analogs of linear 
differential equations with holomorphic or polynomial 
coefficients. Note that in our situation the meaning of a 
polynomial coefficient is not a usual multiplication by a 
polynomial, but the action of a polynomial in the operator $\tau$.

We begin with the regular case and consider an equation (actually, 
a system)
\begin{equation}
dy(t)=P(\tau )y(t)+f(t)
\end{equation}
where for each $z\in \left(\K\right) ^m$, $t\in K$,
\begin{equation}
P(\tau )z=\sum \limits _{k=0}^\infty \pi _kz^{q^k},\quad
f(t)=\sum \limits _{j=0}^\infty \varphi _j\frac{t^{q^j}}{D_j},
\end{equation}
$\pi _k$ are $m\times m$ matrices with elements from $\K$,
$\varphi _j\in \left(\K\right) ^m$, and it is assumed that the series (3.2)
have positive radii of convergence. The action of the operator $\tau $
upon a vector or a matrix is defined component-wise, so that
$z^{q^k}=\left( z_1^{q^k},\ldots ,z_m^{q^k}\right)$ for
$z=(z_1,\ldots ,z_m)$.

We seek a $\F$-linear solution of (3.1) on some neighbourhood of the 
origin, of the form 
\begin{equation}
y(t)=\sum \limits _{i=0}^\infty y_i\frac{t^{q^i}}{D_i},\quad
y_i\in \left(\K\right) ^m,
\end{equation}
where $y_0$ is a given element, so that the ``initial'' condition
for our situation is 
\begin{equation}
\lim \limits _{t\to 0}t^{-1}y(t)=y_0.
\end{equation}

The next theorem, proved in \cite{K00}, is the function field 
analog of the Cauchy theorem from the classical analytic theory of 
differential equations.

\medskip
\begin{teo}
For any $y_0\in \left(\K\right) ^m$ the equation (3.1) has a unique
local solution of the form (3.3), which satisfies (3.4), with the
series having a positive radius of convergence.
\end{teo}

\medskip
Thus, regular equations with Carlitz derivatives behave more or 
less as their classical counterparts. The situation is different 
for singular equations. Let us consider scalar equations of arbitrary order
\begin{equation}
\sum \limits _{j=0}^mA_j(\tau )d^ju=f
\end{equation}
where
$f(t)=\sum \limits _{n=0}^\infty \varphi _n\dfrac{t^{q^n}}{D_n}$,
$A_j(\tau )$ are power series having (as well as the one for $f$)
positive radii of convergence.

We investigate formal solutions of (3.5), of the form
\begin{equation}
u(t)=\sum \limits _{n=0}^\infty u_n\frac{t^{q^n}}{D_n},\quad 
u_n\in \K .
\end{equation}
One can apply an operator series $A(\tau )=\sum \limits _{k=0}^\infty
\alpha _k\tau ^k$ (even without assuming its convergence) to a
formal series (3.6), setting
$$
\tau ^ku(t)=\sum \limits _{n=0}^\infty u_n^{q^k}[n+1]^{q^{k-
1}}\ldots [n+k]\frac{t^{q^{n+k}}}{D_{n+k}},\quad k\ge 1,
$$
and
$$
A(\tau )u(t)=\sum \limits _{l=0}^\infty \frac{t^{q^l}}{D_l}\sum
\limits _{n+k=l}\alpha _ku_n^{q^k}[n+1]^{q^{k-1}}\ldots [n+k]
$$
where the factor $[n+1]^{q^{k-1}}\ldots [n+k]$ is omitted for
$k=0$. These formal manipulations are based on the identity
$$
\tau \left( \frac{t^{q^{i-1}}}{D_{i-1}}\right) 
=[i]\frac{t^{q^i}}{D_i}.
$$
Using also the relation
$$
d\left( \frac{t^{q^i}}{D_i}\right) = \frac{t^{q^{i-1}}}{D_{i-1}},
$$
now we can give a meaning to the notion of a formal solution of 
the equation (3.5).

\medskip
\begin{teo}
Let $u(t)$ be a formal solution (3.6) of the equation (3.5), where
the series for $A_j(\tau )z$, $z\in \K$, and $f(t)$, have
positive radii of convergence. Then the series (3.6) has a positive 
radius of convergence.
\end{teo}

\medskip
This result (proved in \cite{K00}) is in a strong contrast to the 
classical theory. Note that in the $p$-adic case a similar
phenomenon takes place for equations satisfying certain strong conditions
upon zeros of indicial polynomials \cite{Bald,Cla,Put1,Set}. In our case
such a behavior is proved for any equation, which resembles the
(much simpler) case \cite{Put1} of differential equations over a
field of characteristics zero, whose residue field also has 
characteristic zero.

{\bf 3.2.} The equations (3.1) and (3.5) behave like linear 
equations, though they are actually only $\F$-linear. Theorem 3.1 
can be extended \cite{Knl} to the case of strongly nonlinear 
equations (containing self-compositions $y\circ y\circ \cdots 
\circ y$).

On the other hand, it is natural to consider some equations of 
this kind in wider classes of $\F$-linear functions resembling 
meromorphic functions of a complex variables. The set $\mathcal 
R_K$ of locally convergent $\F$-linear holomorphic functions forms 
a non-commutative ring with respect to the composition operation 
(the pointwise multiplication violates the $\F$-linearity). The 
non-commutativity of $\mathcal R_K$ makes the algebraic structures 
related to Carlitz differential equations much more complicated 
compared to their classical counterparts. So far their 
understanding is only at its initial stage. It is known, however, 
that $\mathcal R_K$ can be imbedded into a skew field of 
$\F$-linear ``meromorphic'' series containing terms like 
$t^{q^{-k}}$ (see \cite{Knl}). A deep investigation 
of bi-infinite series of this kind convergent on the whole of $\K$ 
has been carried out by Poonen \cite{Po}.

A specific class of equations with solutions meromorphic in the 
above sense is the class of scalar Riccati-type equations
\begin{equation}
dy(t)=\lambda (y\circ y)(t)+(P(\tau )y)(t)+R(t)
\end{equation}
where $\lambda \in \K$, 
$$
(P(\tau )y)(t)=\sum\limits_{k=1}^\infty p_ky^{q^k}(t),\quad 
R(t)=\sum\limits_{k=0}^\infty r_kt^{q^k},
$$
$p_k,r_k\in \K$ (note that the right-hand side of (3.7) does not 
contain the linear term). The following theorem is proved in 
\cite{Knl}.

\medskip
\begin{teo}
If $0<|\lambda |\le q^{-1/q^2}$, $|p_k|\le q^{-1/q^2}$, $|r_k|\le q^{-1/q^2}$ for 
all $k$, then the equation (3.7) possesses solutions 
of the form
$$
y(t)=ct^{1/q}+\sum\limits_{n=0}^\infty a_nt^{q^n},\quad c,a_n\in 
\K,\ c\ne 0,
$$
where the series converges on the open unit disk $|t|<1$.
\end{teo}

\medskip
\section{Regular Singularity}

{\bf 4.1.} In analysis over $\mathbb C$, a typical class of systems 
with regular singularity at the origin
$\zeta =0$ over $\mathbb C$ consists of systems of the form
\begin{equation}
\zeta y'(\zeta )=\left( B+\sum\limits_{k=1}^\infty
A_k\zeta^k\right) y(\zeta )
\end{equation}
where $B,A_j$ are constant matrices, and the series converges on
a neighbourhood of the origin. Such a system possesses a
fundamental matrix solution of the form $W(\zeta )\zeta^C$ where
$W(\zeta )$ is holomorphic on a neighbourhood of zero, $C$ is a
constant matrix, $\zeta^C=\exp (C\log \zeta )$ is defined by the
obvious power series. Under some additional assumptions regarding
the eigenvalues of the matrix $B$, one can take $C=B$. For
similar results over $\mathbb C_p$ see \cite{DGS}.

In order to investigate such a class of equations in the
framework of $\F$-linear analysis over $K$, one has to go beyond
the class of locally analytic functions. Instead of
power series expansions we can use the expansions in Carlitz
polynomials on the compact ring $O\subset K$. The property of 
local analyticity, if it takes place, can be recovered with the 
use of Theorem 2.3. Note that our approach would fail if we consider
equations over $\K$ instead of $K$ (our solutions may take their
values from $\K$, but they are defined over subsets of $K$). In
this sense our techniques are different from the ones developed
for both the characteristic zero cases.

We begin with the simplest model scalar equation
\begin{equation}
\tau du=\lambda u,\quad \lambda \in \K,
\end{equation}
whose solution may be seen as a function field counterpart of the
power function $t\mapsto t^\lambda$.

We look for a continuous $\F$-linear solution $u(t,\lambda )$ of 
the equation (4.2), with the ``initial condition'' $u(1,\lambda 
)=1$, in the form
\begin{equation}
u(t)=\sum\limits_{i=0}^\infty c_if_i(t),\quad t\in O,
\end{equation}
where $c_0=1$.

It is easy to see that the equation (4.2) has no continuous 
solutions if $|\lambda |\ge 1$. If $|\lambda |<1$, then the 
solution $u(t,\lambda )$ is unique, continuous on $O$, and the 
coefficients from (4.3) have the form
$$
c_n=\prod\limits_{j=0}^{n-1}(\lambda -[j]).
$$

The function $u(t,\lambda )$ is
analytic on $O$ if and only if $\lambda =[j]$ for some $j\ge 0$;
in this case $u(t,\lambda )=u(t,[j])=t^{q^j}$. If $\lambda \ne
[j]$ for any integer $j\ge 0$, then $u(t,\lambda )$ is locally
analytic on $O$ if and only if $\lambda =-x$, and in that case
$u(t,-x)=0$ for $|t|\le q^{-1}$. The relation
$$
u(t^{q^m},\lambda )=u(t,\lambda^{q^m}+[m]),\quad t\in O,
$$
holds for all $\lambda$, $|\lambda |<1$, and for all
$m=0,1,2,\ldots$. For the proofs see \cite{K03}.

Similarly, if in (4.1) $\lambda =(\lambda_{ij})$ is is a 
$m\times m$ matrix with elements from $\K$, and we look for a 
matrix-valued solution of (4.1), then such a solution is given by 
the series (4.3) with the matrix coefficients
$$
c_i=\left\{ \prod\limits_{j=0}^{i-1}(\lambda -[j]I_m)\right\}
c_0,\quad i\ge 1
$$
($I_m$ is the unit matrix), if $|\lambda 
|\stackrel{\text{def}}{=}\max |\lambda_{ij}|<1$.

\bigskip
{\bf 4.2.} The analog, for our situation, of the system (4.1) is 
the system
\begin{equation}
\tau du-P(\tau )u=0
\end{equation}
where $P(\tau )$ is a matrix-valued analytic function, so that 
$P(\tau )z=\sum \limits _{k=0}^\infty \pi _kz^{q^k}$. We assume that 
$|\pi_k|\le \gamma$, $\gamma >0$,
for all $k$, $|\pi_0|<1$. Denote by $g(t)$ a solution of the
equation $\tau dg=\pi_0g$. Let $\lambda_1,\ldots ,\lambda_m\in
\K$ be the eigenvalues of the matrix $\pi_0$.

\begin{teo}
If 
\begin{equation}
\lambda_i-\lambda_j^{q^k}\ne [k],\quad i,j=1,\ldots ,m;\
k=1,2,\ldots,
\end{equation}
then the system (4.4) has a matrix solution
$$
u(t)=W(g(t)),\quad W(s)=\sum\limits_{k=0}^\infty w_ks^{q^k},\quad
w_0=I_m,
$$
where the series for $W$ has a positive radius of convergence.
\end{teo}

\medskip
The paper \cite{K03} contains, apart from the proof of Theorem 
4.1, a discussion of some situations (the Euler type equations) 
where its conditions are violated, as well as of the meaning of 
the conditions (4.5). Here we only mention that in the scalar case 
$m=1$ the condition (4.5) is equivalent to the assumption 
$\pi_0\ne -x$, so that it excludes the case where solutions of the 
equation $\tau dg=\pi_0g$ has pathological properties.

\bigskip
{\bf 4.3.} For the above equation, continuous solutions were found 
as Fourier-Carlitz expansions
\begin{equation}
u(t)=\sum\limits_{n=0}^\infty c_nf_n(t),
\end{equation}
and we had to impose certain conditions upon coefficients of the 
equation, in order to guarantee the uniform convergence of the 
series on $O$ (which is equivalent to the fact that $c_n\to 0$). 
However formally we could write the series (4.6) for the solutions 
without those conditions. Thus, it is natural to ask whether the 
corresponding series (4.6) converge at some points $t\in O$. Note 
that (4.6) always makes sense for $t\in \F [x]$ (for each such $t$ 
only a finite number of terms is different from zero). The 
question is whether the series converges on a wider set; if the 
answer is negative, such a formal solution is called {\it strongly 
singular}.

The available results regarding strong singularity of solutions of 
some equations are based on the following general fact \cite{K03}.

\medskip
\begin{teo}
If $|c_i|\ge \rho >0$ for all $i\ge i_0$ (where $i_0$ is some
natural number), then the function (4.6) is strongly singular.
\end{teo}

\medskip
It follows from Theorem 4.2 that non-trivial formal solutions of
the equation (4.2) with $|\lambda |\ge 1$ are strongly singular. A
more complicated example is provided by the equation
\begin{equation}
(\Delta -[-a])(\Delta -[-b])u=d\Delta u,\quad a,b\in \mathbb Z,
\end{equation}
for Thakur's hypergeometric function ${}_2F_1(a,b;1;t)$.

A holomorphic solution of (4.7) is given by an appropriate 
specialization of (1.5). Classically (over $\mathbb C$), there 
exists the second solution with a logarithmic singularity. Here 
the situation is different. Looking for a solution of the form 
(4.6) we obtain a recursive relation
\begin{multline}
\left( c_{i+2}^{1/q}-c_{i+2}\right) +c_{i+1}^{1/q}[i+1]^{1/q}-
c_{i+1}([i]+[i+1]-[-a]-[-b])\\
-c_i([i]-[-a])([i]-[-b])=0,\quad i=0,1,2,\ldots .
\end{multline}

Taking arbitrary initial coefficients $c_0,c_1\in \K$ we obtain a
solution $u$ defined on $\F [x]$. On each step we have to solve
the equation
\begin{equation}
z^{1/q}-z=v.
\end{equation}
If $|c_i|\le 1$ and $|c_{i+1}|\le 1$, then in the equation for
$c_{i+2}$ we have $|v|<1$.

It can be shown \cite{K03} that the equation (4.9) has a unique 
solution $z_0\in \K$,
for which $|z_0|\le |v|$, and $q-1$ other solutions $z$, $|z|=1$.
It is natural to call
a solution {\it generic} if, starting from a certain step of
finding the coefficients $c_n$, we always take the most frequent
option corresponding to a solution of (4.9) with $|z|=1$. Now 
Theorem 4.2 implies the following fact.

\medskip
\begin{teo}
A generic solution of the equation (4.7) is strongly singular.
\end{teo} 

\medskip
Of course, in some special cases the recursion (4.8) can lead to
more regular solutions, in particular, to the holomorphic
solutions found by Thakur.

\medskip
\section{Polylogarithms and a Zeta Function}

{\bf 5.1.} The Carlitz differential equations can be used for 
defining new special functions with interesting properties. Some 
examples are given in this section.

An analog of the function $-\log (1-t)$ is defined via the equation
\begin{equation}
(1-\tau )du(t)=t,\quad t\in K_\pi ,
\end{equation}
a counterpart of the classical equation $(1-t)u'(t)=1$. The next 
results are taken from \cite{Kzeta} where the equation (5.1) is 
considered for an arbitrary finite place of $\F (x)$.

Let $l_1(t)$ be a $\F$-linear holomorphic solution of (5.1) with 
the zero initial condition (in the sense of (3.4)). Then it is 
easy to show that
\begin{equation}
l_1(t)=\sum\limits_{n=1}^\infty \frac{t^{q^n}}{[n]},
\end{equation}
and the series in (5.2) converges for $|t|\le q^{-1}$.

Note that $l_1(t)$ is different from the well-known Carlitz 
logarithm $\log_C$ (see (1.3)). Analogies motivating the introduction 
of special functions are not so unambiguous, and, for instance, 
from the composition ring viewpoint, $\log_C$ is an analog of 
$e^{-t}$, though in other respects it is a valuable analog of the 
logarithm. By the way, another possible analog of the 
logarithm is a continuous function $u(t)$, $|t|_\pi \le 1$, 
satisfying the equation $\Delta u(t)=t$ (an analog of $tu'(t)=1$) 
and the condition $u(1)=0$. In fact, $u=\mathcal D_1$, the first 
hyperdifferential operator (the definition of $\mathcal D_1$
is given in Sect. 5.2 below); see \cite{Jeong1}.

Now we consider continuous non-holomorphic extensions of $l_1$.

\medskip
\begin{teo}
The equation (5.1) has exactly $q$ continuous solutions on 
$O$ coinciding with (5.2) as $|t|\le q^{-1}$. These 
solutions have the expansions in the Carlitz polynomials
$u=\sum\limits_{i=0}^\infty c_if_i$
where $c_1$ is an arbitrary solution of the equation
$c_1^q-c_1+1=0$, higher coefficients are found from the relation
$$
c_n=\sum\limits_{j=0}^\infty \left( 
c_{n-1}[n-1]\right)^{q^{j+1}},\quad n\ge 2,
$$
and the coefficient $c_0$ is determined by the relation
$$
c_0=\sum\limits_{i=1}^\infty (-1)^{i+1}\frac{c_i}{L_i},
$$
\end{teo}

\medskip
Below we denote by $l_1$ an arbitrary fixed ``branch'' of 
extensions of (5.2).

The polylogarithms $l_n(t)$ are defined recursively by the 
equations
\begin{equation}
\Delta l_n=l_{n-1},\quad n\ge 2,
\end{equation}
which agree with the classical ones $tl_n'(t)=l_{n-1}(t)$.  
Analytic $\F$-linear solutions of (5.3), such that 
$t^{-1}l_n(t)\to 0$ as $t\to 0$, are found easily by induction: 
\begin{equation}
l_n(t)=\sum\limits_{j=1}^\infty \frac{t^{q^j}}{[j]^n},\quad |t|\le 
q^{-1}.
\end{equation}

\medskip
\begin{teo}
For each $n\ge 2$, there exists a unique continuous $\F$-linear 
solution of the equation (5.3) coinciding for $|t|\le 
q^{-1}$ with the polylogarithm (5.4). The solution is given by 
the Carlitz expansion $l_n=\sum\limits_{i=0}^\infty c_i^{(n)}f_i$ with
$$
\left| c_i^{(n)}\right| \le C_nq^{-q^{i-1}},\quad C_n>0,\ 
i\ge 1 ,
$$
$$
c_0^{(n)}=\sum\limits_{i=1}^\infty (-1)^{i+1}\frac{c_i^{(n)}}{L_i},
$$
\end{teo}

\medskip
{\bf 5.2.} Now that the above polylogarithms have been extended 
onto the disk $\{|t|\le 1\}$, we can interpret their values at 
$t=1$ as ``special values'' of a kind of a zeta function. In order 
to define the latter, we introduce the operator $\Delta^{(\alpha )}$, 
$\alpha \in O$, a function field analog of the Hadamard 
fractional derivative $\left( t\frac{d}{dt}\right)^\alpha$ from 
real analysis (see \cite{SKM}).

Denote by $\D_k(t)$, $k\ge 0$, $t\in O$, the sequence of 
hyperdifferentiations defined initially on monomials by the 
relations $\D_0(x^n)=x^n$, $\D_k(1)=0$, $k\ge 1$,
$$
\D_k(x^n)=\binom{n}{k}x^{n-k},
$$
where it is assumed that $\binom{n}{k}=0$ for $k>n$. $\D_k$ is 
extended onto $\F [x]$ by $\F$-linearity, and then onto $O$ by 
continuity \cite{V}. The sequence $\{ \D_k\}$ is an orthonormal 
basis of the space of continuous $\F$-linear functions on $O$ 
\cite{Jeong1,Con}.

Let $\alpha \in O$, $\alpha =\sum\limits_{n=0}^\infty 
\alpha_nx^n$, $\alpha_n\in \F$. Denote $\widehat \alpha =
\sum\limits_{n=0}^\infty (-1)^n\alpha_nx^n$. For an 
arbitrary continuous $\F$-linear function $u$ on $O$ we define 
its ``fractional derivative'' $\Delta^{(\alpha )}u$ at a point 
$t\in O$ by the formula
$$
\left( \Delta^{(\alpha )}u\right) (t)=\sum\limits_{k=0}^\infty 
(-1)^k\D_k(\widehat \alpha )u(x^kt).
$$
The function $\alpha \mapsto \left( \Delta^{(\alpha )}u\right) (t)$ is 
continuous and $\F$-linear. As a function of $t$, $\Delta^{(\alpha 
)}u$ is continuous if, for example, $u$ is H\"older continuous.

Our understanding of $\Delta^{(\alpha )}$ as a kind of a 
fractional derivative is justified by the following properties:
$$
\Delta^{(x^n)}=\Delta^n,\quad n=1,2,\ldots ;
$$
$$
\Delta^{(\alpha )}\left( \Delta^{(\beta )}u\right) (t)= 
\left( \Delta^{(\alpha \beta )}u\right) (t),
$$
for any $\alpha ,\beta \in O$.

\medskip
{\bf 5.3.} We define $\zeta (t)$, $t\in K$, setting $\zeta (0)=0$, 
$$
\zeta (x^{-n})=l_n(1),\quad n=1,2,\ldots ,
$$
and
$$
\zeta (t)=\left( \Delta^{(\theta_0+\theta_1x+\cdots )}l_n\right) 
(1),\quad n=1,2,\ldots ,
$$
if $t=x^{-n}(\theta_0+\theta_1x+\cdots )$, $\theta_j\in \F$. The function 
$\zeta$ is a
continuous $\F$-linear function on $K_x$ with values in $\K$.

In particular, we have 
$$
\zeta (x^m)=\left( \Delta^{m+1}l_1\right) (1),\quad m=0,1,2,\ldots .
$$

The above definition is of course inspired by the classical polylogarithm relation
$$
\left( z\dfrac{d}{dz}\right) \sum\limits_{n=1}^\infty 
\dfrac{z^n}{n^s}=\sum\limits_{n=1}^\infty \dfrac{z^n}{n^{s-1}}.
$$
In contrast to Goss's zeta function defined on natural numbers and 
interpolated onto $\mathbb Z_p$ (see \cite{G2,Th3}), the above 
$\zeta$ is purely an object of the characteristic $p$ arithmetic.

Let us write some relations for special values of our $\zeta$; for 
the details see \cite{Kzeta}.

As we saw,
$$
\sum\limits_{j=1}^\infty \frac{t^{q^j}}{[j]^n}=
\sum\limits_{i=0}^\infty \zeta (x^{-n+i})\D_i(t),\quad |t|\le 
q^{-1}.
$$
Next, let us consider the double sequence $A_{n,r}\in K$, 
$A_{n,1}=(-1)^{n-1}L_{n-1}$,
$$
A_{n,r}=(-1)^{n+r}L_{n-1}\sum\limits_{0<i_1<\ldots 
<i_{r-1}<n}\frac{1}{[i_1][i_2]\ldots [i_{r-1}]},\quad r\ge 2.
$$
These elements appear as the coefficients of the expansion \cite{V} of a 
hyperdifferentiation $\D_r$ in the normalized Carlitz polynomials, 
as well as in the expression \cite{Jeong0} of the operators 
$\Delta^{(n)}$ from (2.1) via the iterations $\Delta^r$. Here we 
have the identity
$$
\zeta (x^{-n})=\sum\limits_{i=1}^\infty (-1)^{i+1}L_i^{-1}
\sum\limits_{r=1}^iA_{i,r}\zeta (x^{r-n})
$$
which may be seen as a distant relative of Riemann's 
functional equation for the classical zeta.

Finally, consider the coefficients $c_i$ of the Carlitz expansion 
of $l_1$ (see Theorem 5.1). They are expressed via zeta values:
$$
c_i=\sum\limits_{r=1}^iA_{i,r}\zeta (x^{r-1}).
$$
By Theorem 5.1, for $i\ge 2$ we have
\begin{equation}
c_i=\sum\limits_{j=0}^\infty (z_i)^{q^j},\quad 
z_i=c_{i-1}^q[i-1]^q.
\end{equation}

The series in (5.5) may be seen as an analog of 
$\sum\limits_jj^{-z}$. This analogy becomes clearer if, for a 
fixed $z\in \K$, $|z|<1$, we consider the set $S$ of all 
convergent power series $\sum\limits_{n=1}^\infty z^{q^{j_n}}$ 
corresponding to sequences $\{ j_n\} \subset \mathbb N$. Let us 
introduce the multiplication $\odot$ in $S$ setting 
$z^{q^i}\odot z^{q^j}=z^{q^{ij}}$ and extending the operation 
distributively (for a similar construction in the framework of 
$q$-analysis in characteristic 0 see \cite{N}). Denoting by 
$\prod\limits_p{}^\odot$ the product in S of elements indexed by 
prime numbers we obtain in a standard way the identity
$$
c_i=\prod\limits_p{}^\odot \sum\limits_{n=0}^\infty (z_i)^{q^{p^n}}
$$
(the infinite product is understood as a limit of the partial 
products in the topology of $\K$), an analog of the Euler 
product formula. It would be interesting to study the algebraic 
structure of $S$ in detail.

\section{Umbral Calculus}

{\bf 6.1.} Classical umbral calculus \cite{RKO,Rom} is a set of algebraic 
tools for obtaining, in a unified way, a rich variety of results 
regarding structure and properties of various polynomial 
sequences. There exists a lot of generalizations extending umbral 
methods to other classes of functions. However there is a 
restriction common to the whole literature on umbral calculus -- 
the underlying field must be of zero characteristic. An attempt to 
mimic the characteristic zero procedures in the positive 
characteristic case \cite{Fe} revealed a number of pathological 
properties of the resulting structures. More importantly, these 
structures were not connected with the existing analysis in 
positive characteristic based on a completely different algebraic 
foundation.

A version of umbral calculus inmplementing such a connection was 
developed by the author \cite{K05}, and we summarize it in this 
section. Its basic notion is motivated by the following identity 
for the non-normalized Carlitz polynomials $e_i=D_if_i$:
\begin{equation}
e_i(st)=\sum\limits_{n=0}^i\binom{i}{n}_Ke_n(t)\{ e_{i-n}(s)\}^{q^n}
\end{equation}
where the ``$K$-binomial coefficients'' $\dbinom{i}{n}_K$ are 
defined as
$$
\binom{i}{n}_K=\frac{D_i}{D_nD_{i-n}^{q^n}}.
$$

Computing the absolute values of the Carlitz factorials directly 
from their definition (1.1), it is easy to show that
$$
\left| \binom{i}{n}_K\right| =1,\quad 0\le n\le i.
$$

In fact, $\dbinom{i}{n}_K\in \F (x)$, and we can consider also 
other places of $\F (x)$, that is other non-equivalent absolute 
values. It can be proved \cite{Khol} that $\dbinom{i}{n}_K$ belongs 
to the ring of integers for any finite place of $\F (x)$.

We see the relation (6.1) as a function field counterpart of the 
classical binomial identity \cite{RKO,Rom} satisfied by many 
classical polynomials. Now, considering a sequence $u_i$ of $\F$-linear
polynomials with coefficients from $\K$, we call it {\it a sequence of 
$K$-binomial type} if $\deg u_i=q^i$ and for all 
$i=0,1,2,\ldots$
\begin{equation}
u_i(st)=\sum\limits_{n=0}^i\binom{i}{n}_Ku_n(t)\left\{ 
u_{i-n}(s)\right\}^{q^n},\quad s,t\in K.
\end{equation} 

As in the conventional umbral calculus, the dual notion is that of 
a delta operator. However, in contrast to the classical situation, 
here the delta operators are only $\F$-linear, not linear.

Denote by $\rho_\lambda$ the operator of multiplicative shift, 
$(\rho_\lambda u)(t)=u(\lambda t)$. We call a linear operator $T$, 
on the $\K$-vector space $\KT$ of all $\F$-linear polynomials, 
{\it invariant} if it commutes with $\rho_\lambda$ for each 
$\lambda \in K$.

A $\F$-linear operator $\delta =\tau^{-1}\delta_0$, where 
$\delta_0$ is a linear invariant operator on $\KT$, is called {\it a 
delta operator} if $\delta_0(t)=0$ and $\delta_0(f)\ne 0$ for $\deg 
f>1$. A sequence $\{ P_n\}_0^\infty$ of $\F$-linear polynomials is 
called {\it a basic sequence} corresponding to a delta operator $\delta 
=\tau^{-1}\delta_0$, if $\deg P_n=q^n$, $P_0(1)=1$, $P_n(1)=0$ for 
$n\ge 1$,
\begin{equation}
\delta P_0=0,\quad \delta P_n=[n]^{1/q}P_{n-1},\ n\ge 1,
\end{equation}
or, equivalently,
\begin{equation}
\delta_0P_0=0,\quad \delta_0P_n=[n]P_{n-1}^q,\ n\ge 1.
\end{equation}

It is clear that $d=\tau^{-1}\Delta$ is a delta operator.
It follows from well-known identities for the Carlitz polynomials 
$e_i$ \cite{G1} (see also (1.10)) that the sequence $\{ e_i\}$ is basic with 
respect to the operator $d$.

\medskip
\begin{teo}
For any delta operator $\delta =\tau^{-1}\delta_0$, there exists a 
unique basic sequence $\{ P_n\}$, which is a sequence of $K$-binomial 
type. Conversely, given a sequence $\{ P_n\}$ of $K$-binomial 
type, define the action of $\delta_0$ on $P_n$ by the relations (6.4), 
extend it onto $\KT$ by linearity and set $\delta =\tau^{-1}\delta_0$. Then
$\delta$ is a delta operator, and $\{ P_n\}$ is the corresponding 
basic sequence.
\end{teo}

\medskip
The analogs of the higher Carlitz difference operators (2.1) in 
the present general context are the operators $\delta_0^{(l)}=\tau^l\delta^l$.
The identity
\begin{equation}
\delta_0^{(l)}P_j=\frac{D_j}{D_{j-l}^{q^l}}P_{j-l}^{q^l}
\end{equation}
holds for any $l\le j$. If $f$ is a $\F$-linear polynomial, $\deg f\le q^n$, 
then a generalized Taylor formula
\begin{equation}
f(st)=\sum\limits_{l=0}^n\frac{\left( \delta_0^{(l)}f\right) 
(s)}{D_l}P_l(t)
\end{equation}
holds for any $s,t\in K$. For the Carlitz polynomials $e_i$, the 
formulas (6.5) and (6.6) are well known \cite{G1}. It is important 
that, in contrast to the classical umbral calculus, the linear 
operators involved in (6.6) are not powers of a single linear 
operator.

Any linear invariant operator $T$ on $\KT$ admits a representation
\begin{equation}
T=\sum\limits_{l=0}^\infty \sigma_l\delta_0^{(l)},\quad 
\sigma_l=\frac{(TP_l)(1)}{D_l}.
\end{equation}
The infinite series in (6.7) becomes actually a finite sum if both 
sides of (6.7) are applied to any $\F$-linear polynomial. 
Conversely, any such series defines a linear invariant operator on 
$\KT$.

Let us consider the case where $\delta =d$, so that 
$\delta_0^{(l)}=\Delta^{(l)}$. The next result leads to new delta 
operators and basic sequences.

\medskip
\begin{teo}
The operator $\theta =\tau^{-1}\theta_0$, where
$$
\theta_0=\sum\limits_{l=1}^\infty \sigma_l\Delta^{(l)},
$$
is a delta operator if and only if
\begin{equation}
S_n\stackrel{\mbox{{\rm \footnotesize def}}}{=}\sum\limits_{l=1}^n
\frac{\sigma_l}{D_{n-l}^{q^l}}\ne 0\ \ \text{for all }\ n=1,2,\ldots .
\end{equation}
\end{teo}

\medskip
{\bf Example 1}.
Let $\sigma_l=1$ for all $l\ge 1$, that is
\begin{equation}
\theta_0=\sum\limits_{l=1}^\infty \Delta^{(l)}.
\end{equation}
Estimates of $|D_n|$ which follow directly from (1.1) show that 
$|S_n|=q^{\frac{q^n-q}{q-1}}$, so that (6.8) is satisfied. 
Comparing (6.9) with a classical formula from \cite{RKO} we may 
see the polynomials $P_n$ for this case as analogs of the Laguerre 
polynomials.

\medskip
{\bf Example 2}. Let $\sigma_l=\dfrac{(-1)^{l+1}}{L_l}$. For this 
case it can be shown \cite{K05} that $S_n=D_n^{-1}$, $n=1,2,\ldots$; 
$\theta_0(t^{q^j})=t^{q^j}$ for all $j\ge 1$ (of course, 
$\theta_0(t)=0$), and $P_0(t)=t$, $P_n(t)=D_n\left( 
t^{q^n}-t^{q^{n-1}}\right)$ for $n\ge 1$.

\medskip
{\bf 6.2.} As in the $p$-adic case \cite{VH,Ver,Rob}, the umbral 
calculus can be used for constructing new orthonormal bases in 
$C_0(O,\K )$.

Let $\{ P_n\}$ be the basic sequence corresponding to a delta 
operator $\delta =\tau^{-1}\delta_0$,
\begin{equation}
\delta_0=\sum\limits_{l=1}^\infty \sigma_l\Delta^{(l)}.
\end{equation}
The sequence $Q_n=\dfrac{P_n}{D_n}$, $n=0,1,2,\ldots$, called the 
normalized basic sequence, satisfies the identity
$$
Q_i(st)=\sum\limits_{n=0}^iQ_n(t)\left\{ Q_{i-n}(s)\right\}^{q^n},
$$
another form of the $K$-binomial property. Though it resembles its 
classical counterpart, the presence of the Frobenius powers is a 
feature specific for the case of a positive characteristic.

\medskip
\begin{teo}
If $|\sigma_1|=1$, $|\sigma_l|\le 1$ for $l\ge 2$, then the 
sequence $\{ Q_n\}_0^\infty$ is an orthonormal basis of the space 
$C_0(O,\K )$ -- for any $f\in C_0(O,\K )$ there is a uniformly 
convergent expansion
$$
f(t)=\sum\limits_{n=0}^\infty \psi_nQ_n(t),\quad t\in O,
$$
where $\psi_n=\left( \delta_0^{(n)}f\right) (1)$, $|\psi_n|\to 0$ 
as $n\to \infty$,
$$
\|f\| =\sup\limits_{n\ge 0}|\psi_n|.
$$
\end{teo}

\medskip
By Theorem 6.3, the Laguerre-type polynomial sequence from 
Example 1 is an orthonormal basis of $C_0(O,\K )$. The sequence 
from Example 2 does not satisfy the conditions of Theorem 6.3.

Note that the conditions of Theorem 6.3 imply that 
$S_n\ne 0$ for all $n$, so that the 
series (6.10) considered in Theorem 6.3 always correspond to delta 
operators.

In \cite{K05} recursive formulas and generating functions for 
normalized basic sequences are also given.

\section{The Weyl-Carlitz Ring and Holonomic Modules}

{\bf 7.1.} The theory of holonomic modules over the Weyl algebra and more 
general algebras of differential or $q$-difference operators is 
becoming increasingly important, both as a crucial part of the 
general theory of D-modules and in view of various applications 
(see, for example, \cite{BK,Cart,S}). Usually, the holonomic 
property of the module corresponding to a system of differential 
equations is a sign of its ``regular'' behavior. Most of the 
classical special functions are associated (see \cite{Cart}) with 
holonomic modules, which helps to investigate their properties.

It is clear from the above 
results that in the positive characteristic case a natural 
counterpart of the Weyl algebra is, for the case of a single 
variable, the ring $\AG_1$ generated by $\tau,d$, and 
scalars from $\K$, with the relations 
\begin{equation}
d\tau -\tau d=[1]^{1/q},\quad \tau \lambda =\lambda^q\tau ,\quad 
d\lambda =\lambda^{1/q}d\ (\lambda \in \K).
\end{equation} 

The ring consists of finite sums
\begin{equation}
a=\sum\limits_{i,j}\lambda_{ij}\tau^id^j,\quad \lambda_{ij}\in \K ,
\end{equation}
and the representation of an element in the form (7.2) is unique.

Basic algebraic properties of $\AG_1$ \cite{K00,B2} are 
similar to those of the Weyl algebra in 
characteristic 0 and quite different from the case of the algebra 
of usual differential operators over a field of positive 
characteristic \cite{Put2}.

The ring $\AG_1$ is left and right Noetherian, without zero 
divisors. $\AG_1$ possesses 
no non-trivial two-sided ideals stable with respect to the mapping
$$
\sum \limits _{i,j}\lambda_{ij}\tau ^id^j\ \mapsto \
\sum \limits _{i,j}\lambda_{ij}^q\tau ^id^j.
$$
The centre of $\AG_1$ is described explicitly in \cite{B2}; it 
contains countably many elements (this corrects an erroneous 
statement from \cite{K00}). In fact, $\AG_1$ belongs to the class 
of generalized Weyl algebras \cite{B1}. A well-developed theory 
available for them enabled Bavula \cite{B2} to classify ideals in 
$\AG_1$, as well as all simple modules over $\AG_1$.

\medskip
A generalization of $\AG_1$ to the case of several variables is 
not straightforward because the Carlitz derivatives $d_s$ and 
$d_t$ do not commute on a monomial $f(s,t)=s^{q^m}t^{q^n}$, if 
$m\ne n$. Moreover, if $m>n$, then $d_s^mf$ is not a polynomial, 
nor even a holomorphic function in $t$ (since the action of $d$ is 
not linear and involves taking the $q$-th root).

A reasonable generalization is inspired by Zeilberger's idea (see 
\cite{Cart}) to study holonomic properties of sequences of 
functions making a transform with respect to the discrete 
variables, which reduces the continuous-discrete case to the 
purely continuous one (simultaneously in all the variables). In 
our situation, if $\{ P_k(s)\}$ is a sequence of $\F$-linear 
polynomials with $\deg P_k\le q^k$, we set
\begin{equation}
f(s,t)=\sum \limits_{k=0}^\infty P_k(s)t^{q^k},
\end{equation}
and $d_s$ is well-defined. In the variable $t$, we consider not 
$d_t$ but the linear operator $\Delta_t$. The latter does not 
commute with $d_s$ either, but satisfies the commutation relations
$$
d_s\Delta_t-\Delta_td_s=[1]^{1/q}d_s,\quad \Delta_t\tau -\tau 
\Delta_t=[1]\tau ,
$$
so that the resulting ring $\AG_2$ resembles a universal 
enveloping algebra of a solvable Lie algebra. 

More generally, denote by $\FF_{n+1}$ the set of all germs of functions of the 
form
\begin{equation}
f(s,t_1,\ldots ,t_n)=\sum\limits_{k_1=0}^\infty \ldots \sum\limits_{k_n=0}^\infty
\sum\limits_{m=0}^{\min (k_1,\ldots ,k_n)}a_{m,k_1,\ldots 
,k_n}s^{q^m}t_1^{q^{k_1}}\ldots t_n^{q^{k_n}}
\end{equation}
where $a_{m,k_1,\ldots ,k_n}\in \K$ are such that all the series 
are convergent on some neighbourhoods of the origin. We do not 
exclude the case $n=0$ where $\FF_1$ will mean the set of all 
$\F$-linear power series $\sum\limits_ma_ms^{q^m}$ convergent on 
a neighbourhood of the origin. $\widehat{\FF}_{n+1}$ will denote 
the set of all polynomials from $\FF_{n+1}$, that is the series 
(7.4) in which only a finite number of coefficients is different 
from zero.

The ring $\AG_{n+1}$ is generated by the operators $\tau 
,d_s,\Delta_{t_1},\ldots \Delta_{t_n}$ on $\FF_{n+1}$, and the 
operators of multiplication by scalars 
from $\K$. To simplify the notation, we write $\Delta_j$ 
instead of $\Delta_{t_j}$ and identify a scalar $\lambda \in \K$ 
with the operator of multiplication by $\lambda$. The operators 
$\Delta_j$ are $\K$-linear, so that
\begin{equation}
\Delta_j\lambda =\lambda \Delta_j,\quad \lambda \in \K ,
\end{equation}
while the operators $\tau ,d_s$ satisfy the commutation relations 
(7.1). In the action of each operator $d_s,\Delta_j$ (acting in a 
single variable), other variables are treated as scalars. The 
operator $\tau$ acts simultaneously on all the variables and 
coefficients. We have the relations involving $\Delta_j$:
\begin{equation}
\Delta_j\tau -\tau \Delta_j=[1]\tau ,\quad 
d_s\Delta_j-\Delta_jd_s=[1]^{1/q}d_s,\quad j=1,\ldots ,n.
\end{equation}

Using the commutation relations (7.1), (7.5), and (7.6), we can 
write any $a\in \AG_{n+1}$, in a unique way, as a finite sum
\begin{equation}
a=\sum c_{l,\mu,i_1,\ldots ,i_n}\tau^ld_s^\mu \Delta_1^{i_1}\ldots 
\Delta_n^{i_n}.
\end{equation} 

\medskip
Let us introduce a filtration in $\AG_{n+1}$ (an analog of the 
Bernstein filtration) denoting 
by $\Gamma_\nu$, $\nu \in \mathbb Z_+$, the $\K$-vector space of 
operators (7.7) with $\max \{ l+\mu +i_1+\cdots +i_n\} \le \nu$ 
where the maximum is taken over all the terms of (7.7). Then 
$\AG_{n+1}$ is a left and right Noetherian filtered ring.

In a standard way (see \cite{Cout}) we define filtered left 
modules over $\AG_{n+1}$. All the basic notions regarding a 
filtered module $M$ (like those of the graded module $\gr (M)$, 
dimension $d(M)$, multiplicity $m(M)$, good filtration etc) are 
introduced just as their counterparts in the theory of modules 
over the Weyl algebra.

If we consider $\AG_{n+1}$ as a left module over itself, then
\begin{equation}
d(\AG_{n+1})=n+2,\quad m(\AG_{n+1})=1.
\end{equation} 
For any finitely generated left $\AG_{n+1}$-module $M$, we have 
$d(M)\le n+2$. By (7.8), this bound cannot be improved in general. 
However, if $I$ is a non-zero left ideal in $\AG_{n+1}$, then
\begin{equation}
d(\AG_{n+1}/I)\le n+1.
\end{equation}
For the module $\widehat{\FF}_{n+1}$ of $\F$-linear polynomials 
(7.4), we have
$$
d\left( \widehat{\FF}_{n+1}\right) =n+1,\quad m\left( \widehat{\FF}_{n+1}\right) 
=n!
$$
The proofs of all these results, as well as the ones given in this 
section below, can be found in \cite{Khol}.

It is natural to call an $\AG_{n+1}$-module $M$ {\it holonomic} if 
$d(M)=n+1$. Thus, $\widehat{\FF}_{n+1}$ is an example of a 
holonomic module.

The next theorem demonstrates, already for the case of 
$\AG_1$-modules, a sharp difference from the case of modules over 
the Weyl algebras. In particular, we see that an analog of the 
Bernstein inequality (see \cite{Cout}) does not hold here without 
some additional assumptions.

\medskip
\begin{teo}
\begin{description}
\item{{\rm (i)}} For any $k=1,2,\ldots$, there exists such a 
nontrivial $\AG_1$-module $M$ that $\dim M=k$ ($\dim$ means the 
dimension over $\K$), that is $d(M)=0$.
\item{{\rm (ii)}} Let $M$ be a finitely generated $\AG_1$-module 
with a good filtration. Suppose that there exists a ``vacuum vector'' $v\in 
M$, such that $d_sv=0$ and $\tau^m(v)\ne 0$ for all 
$m=0,1,2,\ldots$. Then $d(M)\ge 1$.
\end{description}
\end{teo}

\medskip
{\bf 7.3.} Let us consider the case of holonomic submodules of the 
$\AG_{n+1}$-module $\FF_{n+1}$, consisting of $\F$-linear 
functions (7.4) polynomial in $s$ and holomorphic near the origin 
in $t_1,\ldots ,t_n$.

Let $0\ne f\in \FF_{n+1}$,
$$
I_f=\left\{ \varphi \in \AG_{n+1}:\ \varphi (f)=0\right\} .
$$
$I_f$ is a left ideal in $\AG_{n+1}$. The left $\AG_{n+1}$-module 
$M_f=\AG_{n+1}/I_f$ is isomorphic to the submodule 
$\AG_{n+1}f\subset \FF_{n+1}$ -- an element $\varphi (f)\in 
\AG_{n+1}f$ corresponds to the class of $\varphi \in \AG_{n+1}$ in 
$M_f$. A natural good filtration in $M_f$ is induced from that in 
$\AG_{n+1}$. 

As we know (see (7.9)), if $I_f\ne \{ 0\}$, then $d(M_f)\le n+1$. 
We call a function $f$ {\it holonomic} if the module $M_f$ is 
holonomic, that is $d(M_f)=n+1$. The condition $I_f\ne \{ 0\}$ 
means that $f$ is a solution of a non-trivial ``differential equation'' 
$\varphi (f)=0$, $\varphi \in \AG_{n+1}$. The case $n=0$ is quite 
simple.

\medskip
\begin{teo}
If a non-zero function $f\in \FF_1$ satisfies an equation
$\varphi (f)=0$, $0\ne \varphi \in \AG_1$, then $f$ is holonomic.
\end{teo}

\medskip
In particular, any $\F$-linear polynomial of $s$ is holonomic, 
since it is annihilated by $d_s^m$, with a sufficiently large $m$.

If $n>0$, the situation is more complicated. We call 
the module $M_f$ (and the corresponding function $f$) {\it 
degenerate} if $D(M_f)<n+1$ (by the Bernstein inequality, there is 
no degeneracy phenomena for modules over the complex Weyl 
algebra). The simplest example of a degenerate function (for $n=1$)
is $f(s,t_1)=g(st_1)\in \FF_2$ where the function $g$ belongs to
$\FF_1$ and satisfies an equation $\varphi (g)=0$, $\varphi \in 
\AG_1$. It can be shown that $d(M_f)=1$.

In order to exclude the degenerate case, we introduce the notion 
of a non-sparse function.

A function $f\in \FF_{n+1}$ of the form (7.4) is called {\it 
non-sparse} if there exists such a sequence $m_l\to \infty$ that, 
for any $l$, there exist sequences $k_1^{(i)},k_2^{(i)},\ldots 
,k_n^{(i)}\ge m_l$ (depending on $l$), such that $k_\nu^{(i)}\to 
\infty$ as $i\to \infty$ ($\nu =1,\ldots ,n$), and $a_{m,k_1^{(i)},\ldots 
,k_n^{(i)}}\ne 0$.

\medskip
\begin{teo}
If a function $f$ is non-sparse, then $d(M_f)\ge n+1$. If, in 
addition, $f$ satisfies an equation $\varphi (f)=0$, $0\ne \varphi \in 
\AG_{n+1}$, then $f$ is holonomic.
\end{teo}

\medskip
{\bf 7.4} We use Theorem 7.3 to prove that the functions (7.4) 
obtained via the sequence-to-function transform (7.3) or its 
multi-index generalizations, from some well-known sequences of 
polynomials over $K$ are holonomic. In all the cases below the 
non-sparseness is evident, and we have only to prove that the 
corresponding function satisfies a non-trivial Carlitz 
differential equation.

\medskip
a) {\it The Carlitz polynomials}. The transform (7.3) of the 
sequence $\{ f_k\}$ is the Carlitz module function $C_s(t)$; see 
(2.3). It is easy to check that $d_sC_s(t)=C_s(t)$. Therefore 
the Carlitz module function is holonomic, jointly in both its variables.

\medskip
b) {\it Thakur's hypergeometric polynomials}. We consider the 
polynomial case of Thakur's hypergeometric function (1.5), that is
\begin{equation}
{}_lF_\lambda (-a_1,\ldots ,-a_l;-b_1,\ldots ,-b_\lambda 
;z)=\sum\limits_m\frac{(-a_1)_m\ldots (-a_l)_m}{(-b_1)_m\ldots (-b_\lambda 
)_mD_m}z^{q^m}
\end{equation}
where $a_1,\ldots ,a_l,b_1,\ldots ,b_\lambda \in\mathbb Z_+$. It 
is seen from (1.4) that the terms in (7.10), which make sense and 
do not vanish, are those with $m\le \min (a_1,\ldots ,a_l,b_1,\ldots 
,b_\lambda )$. Let the function $f\in \FF_{l+\lambda +1}$ be given 
by
\begin{multline*}
f(s,t_1,\ldots ,t_l,u_1,\ldots ,u_\lambda )\\ =
\sum\limits_{k_1=0}^\infty \ldots \sum\limits_{k_l=0}^\infty
\sum\limits_{\nu_1=0}^\infty \ldots \sum\limits_{\nu_\lambda =0}^\infty
{}_lF_\lambda (-k_1,\ldots ,-k_l;-\nu_1,\ldots ,-\nu_\lambda ;s)\\
\times t_1^{q^{k_1}}\ldots t_l^{q^{k_l}}u_1^{q^{\nu_1}}\ldots 
u_\lambda^{q^{\nu_\lambda}}.
\end{multline*}

It is known (\cite{Th3}, Sect. 6.5) that
\begin{multline}
d_s{}_lF_\lambda (-k_1,\ldots ,-k_l;-\nu_1,\ldots ,-\nu_\lambda 
;s)\\
={}_lF_\lambda (-k_1+1,\ldots ,-k_l+1;-\nu_1+1,\ldots ,-\nu_\lambda +1;s)
\end{multline}
if all the parameters $k_1,\ldots ,k_l,\nu_1,\ldots ,\nu_\lambda$ 
are different from zero. If at least one of them is equal to zero, 
then the left-hand side of (7.11) equals zero. This property implies 
the identity $d_sf=f$, the same as that for the Carlitz module 
function. Thus, $f$ is holonomic.

\medskip
c). {\it $K$-binomial coefficients}. It can be shown \cite{Khol} 
that the $K$-binomial coefficients $\dbinom{k}{m}_K$ (see Sect. 6) 
satisfy the Pascal-type identity
\begin{equation}
\binom{k}{m}_K=\binom{k-1}{m-1}_K^q+\binom{k-1}{m}_K^qD_m^{q-1}
\end{equation}
where $0\le m\le k$ and it is assumed that 
$\dbinom{k}{-1}_K=\dbinom{k-1}{k}_K=0$.

Consider a function $f\in \FF_2$ associated with 
the $K$-binomial coefficients, that is
\begin{equation} 
f(s,t)=\sum\limits_{k=0}^\infty 
\sum\limits_{m=0}^k\binom{k}{m}_Ks^{q^m}t^{q^k}.
\end{equation}
The identity (7.12) implies the equation
$$
d_sf(s,t)=\Delta_tf(s,t)+[1]^{1/q}f(s,t)
$$
for the function (7.13). Therefore $f$ is holonomic.

\end{document}